\documentclass[11pt]{article}

\usepackage[T1]{fontenc}
\usepackage[utf8]{inputenc}
\usepackage{lmodern}
\usepackage[margin=1in]{geometry}
\usepackage{microtype}

\usepackage{parskip}
\usepackage{amsmath,amssymb,amsthm,mathtools,bm}
\newcommand{\Lip}{\mathrm{Lip}}

\newtheorem{theorem}{Theorem}
\newtheorem{lemma}{Lemma}
\newtheorem{proposition}{Proposition}
\newtheorem{corollary}{Corollary}
\newtheorem{definition}{Definition}
\newtheorem{remark}{Remark}

\usepackage{enumitem}

\usepackage{tikz}
\usepackage{pgfplots}
\pgfplotsset{compat=1.18}
\usetikzlibrary{arrows.meta,decorations.pathreplacing,calc}

\usepackage{booktabs}
\usepackage{array}
\usepackage{tabularx}
\usepackage{hyperref}
\hypersetup{colorlinks=true,linkcolor=blue!50!black,citecolor=blue!50!black,urlcolor=blue!50!black}
\usepackage[capitalise]{cleveref}

\usepackage[font=small,labelfont=bf]{caption}
\usepackage{placeins} % for \FloatBarrier if needed

\title{Anchored Implication \& Event-Indexed Fixed Points in Hilbert Spaces: Uniqueness and Quantitative Rates}

\author{
    \begin{tabular}[t]{c}
        \textbf{Faruk Alpay} \\
        \small Department of Logic, Lightcap Institute \\
        \small \texttt{alpay@lightcap.ai}
    \end{tabular}
    \and
    \begin{tabular}[t]{c}
        \textbf{Bugra Kilictas} \\
        \small Computer Engineering, Bahcesehir University \\
        \small \texttt{bugra.kilictas@bahcesehir.edu.tr}
    \end{tabular}
    \and
    \begin{tabular}[t]{c}
        \textbf{Taylan Alpay} \\
        \small Aerospace Engineering, Turkish Aeronautical Association \\
        \small \texttt{s220112602@stu.thk.edu.tr}
    \end{tabular}
}

\date{} % no date

\begin{document}
\maketitle

\begin{abstract}
We develop a synthesis of orthomodular logic (projections as propositions) with operator fixed-point theory in Hilbert spaces. First, we introduce an anchored implication connective $A \Rightarrow^{\mathrm{comm}}_{P} B$, defined semantically so that it is true only when either $A$ is false or else $A$ is true and $B$ is true in a ``commuting'' context specified by a fixed nonzero projection $P$. This connective refines material implication by adding a side condition $[E_B,P]=0$ (commutation of $B$ with the anchor) and reduces to classical implication in the Boolean (commuting) case. Second, we study fixed-point convergence under \emph{event-indexed contractions}. For a \emph{single} nonexpansive (not necessarily linear) map $T$, we prove that the event-indexed condition is \emph{equivalent} to the classical assertion that some power $T^N$ is a strict contraction; thus the ``irregular events'' phrasing does not add generality in that setting. We then present the genuinely more general case of \emph{varying operators} (switching/randomized): if blocks of the evolving composition are contractive with bounded inter-event gaps and a common fixed point exists, we obtain uniqueness and an explicit envelope rate. Finally, with an anchor $P$ that commutes with $T$, the same reasoning ensures convergence on $PH$ under event-indexed contraction on that subspace. We include precise scope conditions, examples, and visual explanations.
\end{abstract}

\section{Introduction}

Hilbert space operator theory and logic based on projections share a rich common ground. Any closed subspace $A$ of a Hilbert space $H$ corresponds to an orthogonal projection $E_A$, which acts as the truth-value operator for the proposition ``the state lies in $A$.'' The lattice of projections is orthomodular, reflecting the non-classical logic of quantum mechanics where not all propositions commute. In classical Boolean logic, implication ($A \implies B$) is defined truth-functionally by $\neg A \vee B$. In an orthomodular lattice, however, there is no single implication connective with all the usual properties; one proposal is the Sasaki implication (also called Sasaki hook) $A \to_S B := A^\perp \vee (A \wedge B)$. The Sasaki arrow coincides with material implication when $A$ and $B$ commute, but in general it lacks some familiar properties. This paper introduces an alternative implication-like connective that incorporates an external anchor projection $P$ to enforce a commutation context. Our anchored implication $A \Rightarrow^{\mathrm{comm}}_{P} B$ is informally read as ``if $A$ then $B$, anchored by $P$,'' and it is designed to hold only in worlds (states) where $B$ lies in the algebra determined by $P$.

On a parallel track, we consider iterative operators on Hilbert space. The classical Banach fixed-point theorem guarantees that a strict contraction mapping on a complete metric space has a unique fixed point and that all iterates converge to it at a geometric rate. However, many operators of interest (e.g., in optimization or quantum dynamics) are nonexpansive---they do not expand distances, but are not contractive at every step---so Banach's theorem does not directly apply.

\section{Preliminaries: Projections, States, and Anchors}

We work in a complex Hilbert space $H$ with inner product $\langle \cdot,\cdot\rangle$ and induced norm $|\cdot|$. Let $B(H)$ denote the algebra of bounded linear operators on $H$. A projection $E = E^* = E^2 \in B(H)$ projects onto some closed subspace of $H$. We write $E_A$ for the orthogonal projection corresponding to a proposition $A$. The range $\mathrm{ran}(E_A) = A$ is the subspace itself, and $E_A^\perp = I - E_A$ is the projection onto the orthogonal complement $A^\perp$. The set of all orthogonal projections on $H$, ordered by $E_A \le E_B$ iff $A \subseteq B$, forms an orthomodular lattice.

A state $\psi \in H$ with $|\psi|=1$ induces a truth valuation $v_\psi$ on propositions by:
\begin{equation}\label{eq:val}
  v_\psi(A) = 
  \begin{cases}
    1,& \text{if } E_A \psi = \psi,\\
    0,& \text{otherwise.}
  \end{cases}
\end{equation}
An \emph{anchor projection} $P$ is a fixed nonzero orthogonal projection on $H$. The commutator $[X,Y] := XY - YX$; note $[E_B,P]=0$ means $E_B$ commutes with $P$.

\paragraph{Lipschitz notation.}
Maps $T:H\to H$ are not assumed linear. We write the global Lipschitz modulus
\[
\Lip(T):=\sup_{x\neq y}\frac{\|Tx-Ty\|}{\|x-y\|}\in[0,\infty].
\]
``Nonexpansive'' means $\Lip(T)\le 1$, and ``$\lambda$-contractive'' means $\Lip(T)\le \lambda<1$. For compositions,
\[
\Lip(S\circ T)\ \le\ \Lip(S)\,\Lip(T),\qquad\text{hence}\quad \Lip(T^p)\le \Lip(T)^p.
\]

\section{Anchored Implication in Orthomodular Logic}

\begin{definition}[Anchored Implication]\label{def:anchored}
Given propositions $A,B$ (projections $E_A,E_B$) and an anchor $P$, define
\begin{equation}\label{eq:anchored-val}
  v_\psi\!\left(A \Rightarrow^{\mathrm{comm}}_{P} B\right) = 1 
  \;\Longleftrightarrow\; 
  \big(v_\psi(A)=0\big)\ \text{or}\ \big(v_\psi(A)=1,\ v_\psi(B)=1,\ [E_B,P]=0\big).
\end{equation}
Thus the implication is vacuously true when $A$ is false; when $A$ is true, it additionally requires $B$ true \emph{and} commutation with $P$ (a global, state-independent condition).
\end{definition}

\begin{proposition}[Reduction to Classical Implication]\label{prop:reduction}
If $P=I$ then $A \Rightarrow^{\mathrm{comm}}_{P} B$ coincides with material implication $A\to B$. More generally, if $[E_A,P]=[E_B,P]=0$ (commuting context), then
\begin{equation}\label{eq:proj-implication}
  E_{\,A \Rightarrow^{\mathrm{comm}}_{P} B} \;=\; E_{\neg A \vee B} \;=\; I - E_A + E_A E_B.
\end{equation}
\emph{Lattice pointer.}
When $E_A$ and $E_B$ commute, the lattice is Boolean: $E_{\neg A\lor B}=E_{\neg A}+E_{A\land B}$ with $E_{\neg A}=I-E_A$ and $E_{A\land B}=E_AE_B$, yielding \eqref{eq:proj-implication}.
\end{proposition}

\begin{remark}[Reduction cases]
\cref{prop:reduction} covers both $P=I$ and the general commuting case $[E_A,P]=[E_B,P]=0$, yielding classical material implication via $E_{\neg A\lor B}=I-E_A+E_AE_B$.
\end{remark}

\begin{remark}[Locality and non-completeness]
Inside the $P$-commutant $\mathcal{C}(P)=\{E\in B(H):[E,P]=0\}$, the anchored implication reduces to the classical connective. This is a \emph{sharp but local} observation: it does not entail, and we do not claim, any global completeness or representation theorem for the full orthomodular logic beyond $\mathcal{C}(P)$.
\end{remark}

\begin{lemma}[Monotonicity in the Antecedent]\label{lem:monotonicity}
If $E_A \le E_{A'}$ and $[E_B,P]=0$, then $A \Rightarrow^{\mathrm{comm}}_{P} B \le A' \Rightarrow^{\mathrm{comm}}_{P} B$.
\begin{proof}
Fix $\psi$. If $v_\psi(A)=0$ then $v_\psi(A \Rightarrow^{\mathrm{comm}}_{P} B)=1$, and $v_\psi(A')\ge v_\psi(A)$ forces $v_\psi(A' \Rightarrow^{\mathrm{comm}}_{P} B)=1$. If $v_\psi(A)=1$ then $v_\psi(A')=1$. Under the \emph{state-independent} side condition $[E_B,P]=0$, \eqref{eq:anchored-val} reduces both implications to $v_\psi(B)$, hence the claim.
\end{proof}
\end{lemma}

\paragraph{Worked example (non-commuting anchor; failure despite $A$ and $B$ true).}
Let $H=\mathbb{R}^2$, $A=B=\mathrm{span}\{e_1\}$, $E_B=\begin{psmallmatrix}1&0\\0&0\end{psmallmatrix}$, and let
$P=\tfrac12\begin{psmallmatrix}1&1\\1&1\end{psmallmatrix}$ (onto $(1,1)/\sqrt2$). Then
\[
[E_B,P]=\tfrac12\begin{psmallmatrix}0&1\\-1&0\end{psmallmatrix}\neq 0.
\]
For $\psi=e_1$, $v_\psi(A)=v_\psi(B)=1$ but $v_\psi(A \Rightarrow^{\mathrm{comm}}_{P} B)=0$ because the commutation side condition fails.

\paragraph{Worked commuting example (explicit matrices).}
Let $E_A=\begin{psmallmatrix}1&0\\0&0\end{psmallmatrix}$ and $E_B=\begin{psmallmatrix}1&0\\0&1\end{psmallmatrix}$ (so $E_A,E_B$ commute) and take $P=\begin{psmallmatrix}1&0\\0&0\end{psmallmatrix}$; then $[E_A,P]=[E_B,P]=0$. We compute
\[
I-E_A+E_AE_B=\begin{psmallmatrix}1&0\\0&1\end{psmallmatrix}-\begin{psmallmatrix}1&0\\0&0\end{psmallmatrix}+\begin{psmallmatrix}1&0\\0&0\end{psmallmatrix}=\begin{psmallmatrix}1&0\\0&1\end{psmallmatrix}=E_{\neg A\lor B}.
\]

\begin{proposition}[No Uniform Classical Synonym]\label{prop:no-synonym}
There is no Boolean $F(A,B)$ such that for all anchors $P$,
$A \Rightarrow^{\mathrm{comm}}_{P} B \equiv F(A,B)$.
\begin{proof}
If $[E_A,P]=[E_B,P]=0$, then by \cref{prop:reduction} one has $A \Rightarrow^{\mathrm{comm}}_{P} B \equiv \neg A \lor B$.
If $[E_B,P]\neq 0$, then \eqref{eq:anchored-val} yields $v_\psi(A \Rightarrow^{\mathrm{comm}}_{P} B)=v_\psi(\neg A)$ for all $\psi$, all $A,B$. The table shows the mismatch:
\[
\begin{array}{c|c|c}
(A,B) & v(A \Rightarrow^{\mathrm{comm}}_{P} B)\ \ \text{if}\ [E_B,P]\neq 0 & v(\neg A\vee B)\\
\hline
00 & 1 & 1\\
01 & 1 & 1\\
10 & 0 & 0\\
11 & 0 & \mathbf{1}\ \text{mismatch}
\end{array}
\]
Hence no single Boolean $F(A,B)$ agrees with both anchor regimes.
\end{proof}
\end{proposition}

\begin{remark}[Authors' note: Prop.~\ref{prop:no-synonym}]
The truth-table argument is definitive. A fully formalized variant for effect-valued/alternate semantics is deferred to \emph{Appendix~B} (thresholded effects), which preserves the commuting reduction and clarifies the scope of this claim.
\end{remark}

\subsection*{Algebraic notes (with proofs)}
Let $\mathcal{C}(P)=\{E\in B(H):[E,P]=0\}$ denote the $P$-commutant.

\begin{lemma}[Partial residuation in $\mathcal{C}(P)$]\label{lem:partial-residuation}
If $A,B,X\in\mathcal{C}(P)$, then
\[
X \le A \Rightarrow^{\mathrm{comm}}_{P} B \quad \Longleftrightarrow \quad A\wedge X \le B .
\]
\begin{proof}
By \cref{prop:reduction}, $E_{A \Rightarrow^{\mathrm{comm}}_{P} B}=I-E_A+E_AE_B$.
For commuting projections, $E_{A\wedge X}=E_AE_X$. Then $E_X\le I-E_A+E_AE_B$ iff $E_AE_X\le E_B$, i.e., $A\wedge X\le B$.
\end{proof}
\end{lemma}

\paragraph{Sequent rules (soundness).}
\[
\frac{\ }{A \Rightarrow^{\mathrm{comm}}_{P} B}\quad \textsc{(AI-Intro-Vac)}\quad\text{if }\neg A,
\qquad
\frac{A\quad B\quad [E_B,P]=0}{A \Rightarrow^{\mathrm{comm}}_{P} B}\quad \textsc{(AI-Intro-Comm)},
\]
\[
\frac{A \Rightarrow^{\mathrm{comm}}_{P} B\quad A\quad [E_B,P]=0}{B}\quad \textsc{(AI-Elim)},
\qquad
\frac{A\le A'\quad [E_B,P]=0}{A \Rightarrow^{\mathrm{comm}}_{P} B\ \le\ A' \Rightarrow^{\mathrm{comm}}_{P} B}\quad \textsc{(AI-Mono(A))}.
\]

\paragraph{Rule illustrations.}
\emph{AI-Intro-Vac:} If $v_\psi(A)=0$, then by \eqref{eq:anchored-val} the implication holds regardless of $B$ or $P$. \emph{AI-Intro-Comm:} If $v_\psi(A)=v_\psi(B)=1$ and $[E_B,P]=0$, then \eqref{eq:anchored-val} yields truth. \emph{AI-Elim:} If $A \Rightarrow^{\mathrm{comm}}_{P} B$ and $A$ hold and $[E_B,P]=0$, \eqref{eq:anchored-val} forces $v_\psi(B)=1$. \emph{AI-Mono(A):} If $E_A\le E_{A'}$ and $[E_B,P]=0$, evaluations agree as in \cref{lem:monotonicity}.

\section{Fixed Points: Single Map vs.\ Varying Operators}\label{sec:fixedpoints}

\begin{definition}[Event-indexed contraction]\label{def:event}
A nonexpansive $T:H\to H$ has an \emph{event-indexed $\lambda$-contraction} if there exist $0<\lambda<1$ and an index $n_1\in\mathbb{N}$ such that
\begin{equation}\label{eq:event-contract}
  \|T^{n_1}x - T^{n_1}y\| \le \lambda\, \|x-y\| \qquad \forall x,y \in H,
\end{equation}
equivalently $\Lip(T^{n_1})\le \lambda$.
\end{definition}

\begin{proposition}[Single-map equivalence]\label{prop:eq-single}
Let $T$ be nonexpansive. The following are equivalent:
\begin{enumerate}[label=(\alph*),itemsep=2pt,leftmargin=1.2em]
\item There exist $\lambda\in(0,1)$ and $n_1\in\mathbb{N}$ such that \eqref{eq:event-contract} holds.
\item There exists $N\in\mathbb{N}$ and $\lambda\in(0,1)$ such that $T^N$ is $\lambda$-contractive on $H$ (i.e., $\Lip(T^N)\le \lambda$).
\end{enumerate}
Moreover, if either holds, then for every $p\ge N$,
\[
\Lip(T^{p})\ \le\ \lambda.
\]
Indeed, write $p=N+k$ and use $\Lip(T^{N+k})\le \Lip(T^N)\,\Lip(T)^k\le \lambda$.
\end{proposition}

\begin{remark}[Authors' note: Prop.~\ref{prop:eq-single} and Banach]
$T$ is not assumed linear; submultiplicativity of $\Lip(\cdot)$ suffices. As a corollary, Banach’s fixed-point theorem applies to $T^N$, giving a unique fixed point for $T^N$ which, by \cref{lem:inheritance}, is the unique fixed point of $T$.
\end{remark}

\begin{lemma}[Fixed point inheritance]\label{lem:inheritance}
If $T^N$ is a contraction on $H$, it has a unique fixed point $z$, which is also the unique fixed point of $T$.
\begin{proof}
$T^{N}(Tz)=T(T^{N}z)=Tz$, so $Tz$ is fixed by $T^N$. By uniqueness, $Tz=z$, hence $z$ is fixed by $T$. If $w$ is any fixed point of $T$, then $T^Nw=w$, so $w=z$ by uniqueness.
\end{proof}
\end{lemma}

\begin{corollary}[Classical rate]\label{cor:classical}
Assume \cref{prop:eq-single}. Then $T$ has a unique fixed point $z$ and, for all $n\ge N$,
\[
\|T^n x - z\| \le \lambda^{\,n-N+1}\,\|T^{N-1}x - z\|.
\]
\end{corollary}

\subsection{Varying operators: where event-indexing adds generality}

Let $(T_t)_{t\ge 1}$ be a sequence of nonexpansive maps on $H$, and write forward products
\[
\Phi(n{:}m):=T_n T_{n-1}\cdots T_m,\qquad n\ge m.
\]
Assume there exists a nonempty set $\mathrm{Fix}:=\bigcap_{t\ge1}\mathrm{Fix}(T_t)$ (a common fixed point set).

\begin{theorem}[Event-indexed block contractions]\label{thm:main}
Suppose there exist $\lambda\in(0,1)$, $M\in\mathbb{N}$, and indices $1\le n_1<n_2<\cdots$ with $n_{k+1}-n_k\le M$ such that each \emph{block} is contractive:
\begin{equation}\label{eq:block-contract}
\|\Phi(n_k{:}n_{k-1}{+}1) u - \Phi(n_k{:}n_{k-1}{+}1) v\|\ \le\ \lambda\,\|u-v\|\quad\forall u,v\in H,\ \forall k,
\end{equation}
with $n_0:=0$ (equivalently, $\Lip\!\big(\Phi(n_k{:}n_{k-1}{+}1)\big)\le \lambda$). Then for any $x_0\in H$ the orbit $x_n:=\Phi(n{:}1)x_0$ converges to the \emph{unique} $z\in\mathrm{Fix}$, and for all $n\ge n_1$,
\begin{equation}\label{eq:rate-main}
\|x_n - z\|\ \le\ \lambda^{\,1 + \left\lfloor \frac{n - n_1}{M} \right\rfloor}\, \|x_{n_1}-z\|.
\end{equation}
\begin{proof}
Pick $z\in\mathrm{Fix}$. For any $t$, $T_t z=z$, so between events
$\|x_{t+1}-z\|=\|T_{t+1}x_t - T_{t+1}z\|\le \|x_t-z\|$.
At an event block, \eqref{eq:block-contract} yields
$\|x_{n_k}-z\|\le \lambda\,\|x_{n_{k-1}}-z\|$, hence $\|x_{n_k}-z\|\le \lambda^k\|x_0-z\|$ by induction. For $n\in[n_k,n_{k+1})$, monotonicity gives $\|x_n-z\|\le \|x_{n_k}-z\|\le \lambda^k\|x_0-z\|$. Since $n_{k+1}-n_k\le M$, $k\ge 1+\lfloor (n-n_1)/M\rfloor$ for $n\ge n_1$, establishing \eqref{eq:rate-main}. If $w\in\mathrm{Fix}$, then for each $k$,
$\|w-z\|\le \lambda\,\|w-z\|$, forcing $w=z$.
\end{proof}
\end{theorem}

\begin{lemma}[Common fixed points: sufficient conditions]\label{lem:fix-nonempty}
The set $\mathrm{Fix}=\bigcap_{t\ge1}\mathrm{Fix}(T_t)$ is nonempty under any of the following settings:
\begin{enumerate}[label=(\roman*),leftmargin=1.2em]
\item $T_t=P_{C_t}$ are metric projections onto nonempty closed convex sets $C_t$ with $\bigcap_t C_t\neq\varnothing$;
\item $T_t=\mathrm{prox}_{\gamma_t f_t}$ with a shared minimizer $x^\star\in\bigcap_t \arg\min f_t$;
\item $T_t=J_{\gamma_t A_t}$ (resolvents of maximally monotone $A_t$) with a common zero $x^\star\in\bigcap_t A_t^{-1}(0)$.
\end{enumerate}
\begin{proof}
In (i)–(iii), any $z$ in the stated intersection satisfies $T_t z=z$ for all $t$, so $z\in\mathrm{Fix}$.
\end{proof}
\end{lemma}

\begin{remark}[Authors' note: Thm.~\ref{thm:main}]
The nonemptiness of $\mathrm{Fix}$ is essential. \cref{lem:fix-nonempty} lists standard verification paths (projections, proximal maps, resolvents).
\end{remark}

\begin{theorem}[Heterogeneous blocks and variable gaps]\label{thm:hetero}
Let $0<\lambda_k<1$ and let $1\le n_1<n_2<\cdots$, with block lengths $N_k:=n_k-n_{k-1}\in\mathbb{N}$ ($n_0:=0$). Assume $\Lip\!\big(\Phi(n_k{:}n_{k-1}{+}1)\big)\le \lambda_k$ for all $k$, each $T_t$ is nonexpansive, and $\mathrm{Fix}\neq\varnothing$. Then for any $x_0$ and $z\in\mathrm{Fix}$, letting $k(n):=\max\{k:\ n_k\le n\}$,
\[
\|x_n-z\|\ \le\ \Big(\prod_{j=1}^{k(n)}\lambda_j\Big)\,\|x_{n_1}-z\|\qquad(n\ge n_1).
\]
Moreover,
\[
\limsup_{n\to\infty}\frac{1}{n}\ln\|x_n-z\|\ \le\ \limsup_{K\to\infty}\frac{\sum_{j=1}^K \ln \lambda_j}{\sum_{j=1}^K N_j}.
\]
If, in addition, $N_k\le M$ for all $k$, then
\[
\|x_n-z\|\ \le\ \Bigg(\prod_{j=1}^{1+\left\lfloor\frac{n-n_1}{M}\right\rfloor}\lambda_j\Bigg)\,\|x_{n_1}-z\|
\]
and
\[
\limsup_{n\to\infty}\frac{1}{n}\ln\|x_n-z\|\ \le\ \frac{1}{M}\,\limsup_{K\to\infty}\frac{1}{K}\sum_{j=1}^K \ln\lambda_j.
\]
\begin{proof}
Between events, nonexpansiveness yields monotone decrease of $\|x_n-z\|$; each block multiplies the distance by at most $\lambda_k$, giving the product bound. For the slope bound, note $n\ge \sum_{j=1}^{k(n)}N_j$ and apply $\ln$ to the product. The $M$-bounded case uses $\sum_{j=1}^{K}N_j\le MK$.
\end{proof}
\end{theorem}

\begin{remark}[Examples illustrating sharpness]\label{rem:hetero-examples}
\emph{(i) Alternating factors, varying gaps.)} Let $\lambda_{2j-1}=0.7$, $\lambda_{2j}=0.9$, with block lengths $N_{2j-1}=2$, $N_{2j}=3$. Then after $K$ blocks,
$\|x_{n_K}-z\|\le (0.7\cdot 0.9)^{K/2}\|x_{n_0}-z\|$ and
\[
\frac{1}{n_K}\ln\|x_{n_K}-z\| \le \frac{\frac{K}{2}(\ln 0.7+\ln 0.9)}{\sum_{j\le K}N_j}
\to \frac{\ln 0.7+\ln 0.9}{2+3}<0.
\]
\emph{(ii) Borderline slope $0$.} Take $\lambda_k=1-\frac{1}{k^2}$ and $N_k=k^2$. Then $\sum_{k}\ln\lambda_k$ converges while $\sum_k N_k$ diverges, so the log-slope bound tends to $0^{-}$, showing the envelope can be as flat as desired without becoming positive.
\end{remark}

\begin{definition}[Event schedule and log-envelope]\label{def:envelope}
Given $\lambda\in(0,1)$, a gap bound $M\in\mathbb{N}$ and first event time $n_1$, define the \emph{log-envelope}
\[
E(n)\ :=\ \lambda^{\,1 + \left\lfloor \frac{n - n_1}{M} \right\rfloor}\!,\qquad n\ge n_1,
\]
whose logarithmic slope is $\frac{\ln\lambda}{M}$ between events.
\end{definition}

\begin{lemma}[Envelope tightness for periodic events]\label{lem:tight}
If events occur periodically with period $M$ and during each event the map is an exact $\lambda$-contraction while all inter-event steps are isometries (i.e., $\Lip(\cdot)=1$), then there exist $x_0,z$ with $z$ fixed such that
\[
\|x_n-z\|\ =\ E(n)\,\|x_{n_1}-z\|\quad\text{for all }n\ge n_1.
\]
In particular, the envelope slope $\ln(\lambda)/M$ is unimprovable.
\begin{proof}
Let inter-event steps be isometries and the event step be $\lambda I$. Then $\|x_{n_k}-z\|=\lambda^k\|x_{n_0}-z\|$ and distances are constant between events; this realizes $E(n)$ with equality.
\end{proof}
\end{lemma}

\paragraph{Non-periodic heterogeneous schedule realizing the product envelope (exact).}
Let $H=\mathbb{R}^2$ and fix any sequences $(\lambda_k)\subset(0,1)$ and $(N_k)\subset\mathbb{N}$. For each block $k$, define $N_k-1$ inter-event steps as isometries (e.g., rotations $R_{\pi/2}$), and the last step as $\lambda_k I$. Then for $n\in[n_k,n_{k+1})$,
\[
\|x_n-z\|\ =\ \Big(\prod_{j=1}^{k}\lambda_j\Big)\,\|x_{n_0}-z\|,
\]
i.e., the heterogeneous product envelope is achieved with equality at all times.

\paragraph{Concrete periodic example (reproducible iterates).}
Let $H=\mathbb{R}^2$, period $M=4$ with
$T_1=T_2=T_3=R_{\pi/2}$ and $T_4=\alpha I$ ($\alpha=0.8$), repeated. The block map is $\Phi(4{:}1)=\alpha\, R_{\pi/2}^3$. With $z=0$ and $x_0=(1,0)$,
\[
x_{4k+r}\ =\ R_{\pi/2}^{\,r}\,\big(\alpha^k R_{\pi/2}^{\,3k}x_0\big),\qquad
\|x_{4k+r}-z\|\ =\ \alpha^k,\quad r\in\{0,1,2,3\}.
\]

\begin{figure}[htbp]
\centering
\begin{tikzpicture}
\begin{axis}[
  width=\textwidth,
  height=0.36\textwidth,
  xmin=0, xmax=32,
  ymin=0, ymax=1.05,
  title={Effect of periodic contraction blocks on convergence},
  xlabel={Iteration $n$ (steps)},
  ylabel={Distance $\|x_n - z\|$ (norm units)},
  grid=both,
  major grid style={gray!35},
  minor grid style={gray!20},
  tick align=outside,
  legend style={draw=none, fill=white, font=\small, at={(0.02,0.18)}, anchor=west},
  clip=false
]
  % Periodic events at n = 4,8,12,16,20,24,28
  \addplot+[const plot, very thick, mark=square*]
    coordinates {
      (0,1.0) (1,1.0) (2,1.0) (3,1.0) (4,0.8)
      (5,0.8) (6,0.8) (7,0.8) (8,0.64)
      (9,0.64) (10,0.64) (11,0.64) (12,0.512)
      (13,0.512) (14,0.512) (15,0.512) (16,0.4096)
      (17,0.4096) (18,0.4096) (19,0.4096) (20,0.32768)
      (21,0.32768) (22,0.32768) (23,0.32768) (24,0.262144)
      (25,0.262144) (26,0.262144) (27,0.262144) (28,0.2097152)
      (29,0.2097152) (30,0.2097152) (31,0.2097152)
    };
  \addlegendentry{with event-indexed blocks ($\lambda=0.8$, $M=4$; $x_0=(1,0)$, $z=0$)};
  \addplot+[domain=0:32, samples=2, dashed, very thick, mark=*] {1};
  \addlegendentry{no contractive blocks (nonexpansive only)};
\end{axis}
\end{tikzpicture}
\caption{Distance vs.\ iteration in norm units for the periodic example with $x_0=(1,0)$ and $z=0$. Contractive blocks (solid squares at $n=4,8,12,\dots$) cause geometric drops; without blocks (dashed circles), no decay is guaranteed.}
\label{fig:fig1}
\end{figure}
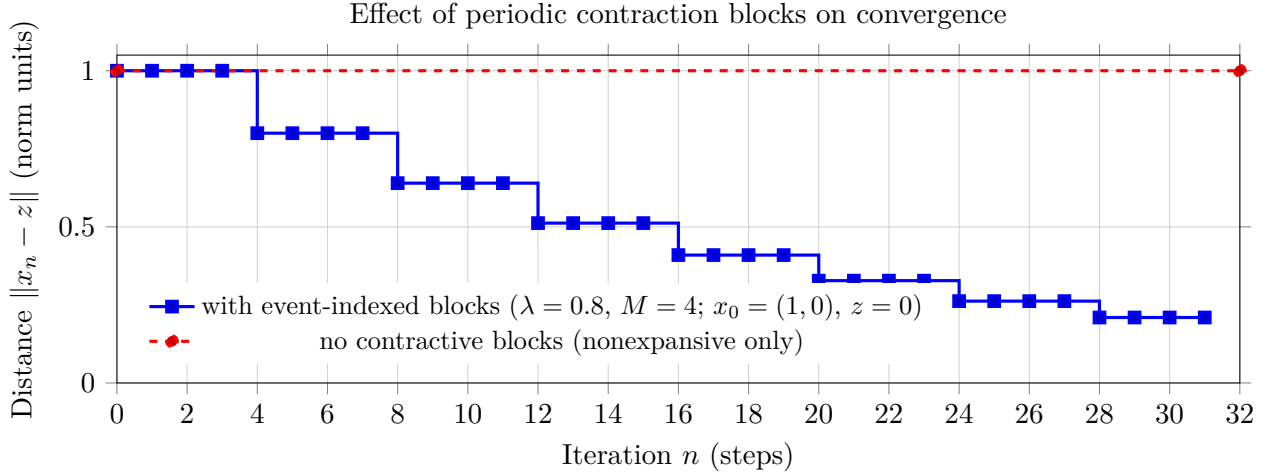

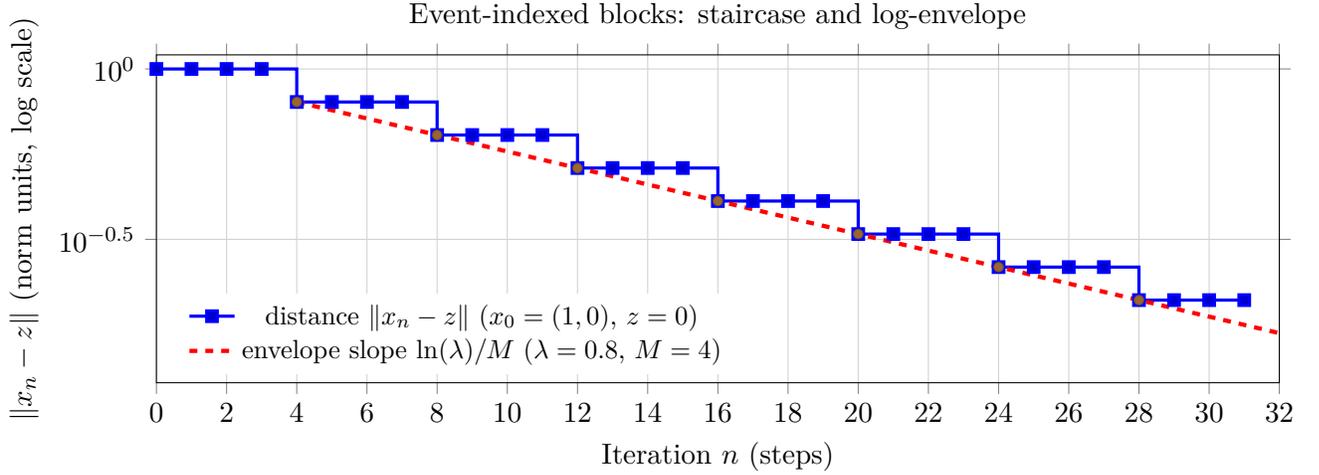
\begin{figure}[htbp]
\centering
\begin{tikzpicture}
\begin{axis}[
  width=\textwidth,
  height=0.36\textwidth,
  xmin=0, xmax=32,
  ymin=0.12, ymax=1.1,
  ymode=log,
  title={Event-indexed blocks: staircase and log-envelope},
  xlabel={Iteration $n$ (steps)},
  ylabel={$\|x_n - z\|$ (norm units, log scale)},
  grid=both,
  major grid style={gray!35},
  minor grid style={gray!20},
  tick align=outside,
  legend style={draw=none, fill=white, font=\small, at={(0.02,0.15)}, anchor=west},
  clip=false
]
  \addplot+[const plot, very thick, mark=square*]
    coordinates {
      (0,1.0) (1,1.0) (2,1.0) (3,1.0) (4,0.8)
      (5,0.8) (6,0.8) (7,0.8) (8,0.64)
      (9,0.64) (10,0.64) (11,0.64) (12,0.512)
      (13,0.512) (14,0.512) (15,0.512) (16,0.4096)
      (17,0.4096) (18,0.4096) (19,0.4096) (20,0.32768)
      (21,0.32768) (22,0.32768) (23,0.32768) (24,0.262144)
      (25,0.262144) (26,0.262144) (27,0.262144) (28,0.2097152)
      (29,0.2097152) (30,0.2097152) (31,0.2097152)
    };
  \addlegendentry{distance $\|x_n - z\|$ ($x_0=(1,0)$, $z=0$)};
  \addplot+[domain=4:32, samples=200, dashed, ultra thick, mark=none]
    {0.8 * exp(ln(0.8)/4 * (x - 4))};
  \addlegendentry{envelope slope $\ln(\lambda)/M$ ($\lambda=0.8$, $M=4$)};
  \addplot+[only marks, mark=*, mark size=1.8pt]
    coordinates {(4,0.8) (8,0.64) (12,0.512) (16,0.4096) (20,0.32768) (24,0.262144) (28,0.2097152)};
\end{axis}
\end{tikzpicture}
\caption{Each event block reduces distance by $\lambda$; between events, nonexpansiveness keeps $\|x_n-z\|$ non-increasing. The dashed envelope has slope $\ln(\lambda)/M$. The iterates shown are exactly $x_{4k+r}=R_{\pi/2}^{\,r}(\alpha^k R_{\pi/2}^{\,3k}x_0)$ with $x_0=(1,0)$, $z=0$.}
\label{fig:fig2}
\end{figure}

\begin{proposition}[No-events $\Rightarrow$ possible non-convergence]\label{prop:noevents}
If no contractive block ever occurs (i.e., \eqref{eq:block-contract} never holds), convergence can fail.
\begin{proof}[Proof sketch]
Let $T_t\equiv R_{\pi/2}$ (a $90^\circ$ rotation). Then $\Lip(T_t)=1$ and every block is nonexpansive but noncontractive. For $x_0\neq 0$, the orbit $x_n=R_{\pi/2}^n x_0$ cycles on the unit circle and does not converge. Thus, without contractive blocks, convergence need not occur.
\end{proof}
\end{proposition}

\begin{remark}[Milder regularity without explicit block-contractivity]
Even if \eqref{eq:block-contract} fails, convergence may occur under \emph{averagedness} and Fejér-type conditions. For instance, if all $T_t$ are averaged with respect to a nonempty closed convex set $C$ and the iterates are Fejér-monotone to $C$, then $\|x_{n+1}-z\|\le \|x_n-z\|$ for all $z\in C$; with Opial's property one obtains weak convergence of $(x_n)$ to a point in $C$, and stronger conditions (e.g., demiregularity/metric subregularity) yield strong or even linear convergence.
\end{remark}

\begin{remark}[Not an iff]\label{rem:not-iff}
Event-indexed block contractivity is a convenient \emph{sufficient} recipe for uniqueness and rate; no ``if and only if'' is claimed.
\end{remark}

\section{Anchored Fixed Points on Invariant Subspaces}

Suppose $P$ is an anchor with a single $T$ commuting, $TP=PT$. Then $PH$ is $T$-invariant and we can analyze $PT^n x = T^n (P x)$.

\begin{theorem}[Anchored convergence under invariance]\label{thm:anchored}
Let $T$ be nonexpansive and $P\neq 0$ with $TP=PT$. If there exist $0<\lambda<1$ and $N$ such that $(T|_{PH})^N$ is $\lambda$-contractive on $PH$, then $T|_{PH}$ has a unique fixed point $z\in PH$ and
\begin{equation}\label{eq:anchored-rate}
\|PT^n x - z\|\le \lambda^{\,n-N+1}\,\|T^{N-1}(Px) - z\|\qquad(n\ge N).
\end{equation}
\begin{proof}[Full proof]
\emph{Invariance and nonexpansiveness on $PH$.} Since $TP=PT$, for any $u\in H$ we have $T(Pu)=P(Tu)\in PH$, so $PH$ is invariant and $T|_{PH}$ is well-defined and nonexpansive: $\Lip(T|_{PH})\le \Lip(T)\le 1$.

\emph{Contractive power on $PH$.} By hypothesis, $\Lip((T|_{PH})^N)\le \lambda<1$. Banach’s theorem on the complete metric space $(PH,\|\cdot\|)$ yields a unique $z\in PH$ with $T|_{PH}(z)=z$.

\emph{Uniqueness on $PH$.} If $w\in PH$ and $T|_{PH}(w)=w$, then $(T|_{PH})^N(w)=w$, so by contraction of $(T|_{PH})^N$ we have $w=z$.

\emph{Rate bound.} For $n\ge N$,
\[
\|PT^n x - z\|=\|T^n(Px)-z\|=\|(T|_{PH})^n(Px)-z\|\le \lambda^{\,n-N+1}\,\|(T|_{PH})^{N-1}(Px)-z\|,
\]
by the same telescoping estimate as in \cref{cor:classical}.

\emph{Anchoring effect.} Note that $T$ may have multiple global fixed points in $H$, yet $T|_{PH}$ has a \emph{unique} fixed point $z$; for instance, $T(x,y)=(x,0.5y)$ has $\mathrm{Fix}(T)=\{(x,0):x\in\mathbb{R}\}$ whereas for $P=\mathrm{diag}(0,1)$, $T|_{PH}$ has the unique fixed point $(0,0)$.
\end{proof}
\end{theorem}

\noindent\textit{Mini-scope.} Holds when $TP=PT$ and $(T|_{PH})^N$ is contractive; not claimed when $TP\neq PT$.

\section{Related Work}\label{sec:related}

\paragraph{Implication connectives in orthomodular lattices.}
Beyond the Sasaki arrow, a variety of implication-like operations have been studied in orthomodular lattices and quantum logics, including residuated or commutation-conditioned connectives and modal/operational viewpoints. Early logical treatments appear in Piziak’s work on implication algebras in orthomodular settings \cite{piziak1974}, with broader lattice-theoretic background in Kalmbach \cite{kalmbach1983}. Our anchored implication adds a \emph{global} commutation side condition tied to a fixed anchor $P$, ensuring classical behavior inside the $P$-commutant while exposing a principled failure mode when $[E_B,P]\neq 0$ (cf.\ \cref{prop:reduction,prop:no-synonym}). This aligns with operational perspectives where additional structure (e.g., commutativity constraints) governs admissible inference; see also Gudder’s expository discussion on projections-as-propositions in quantum computation \cite{gudder2003}.

\paragraph{Convergence of compositions of nonexpansive/averaged operators.}
Fixed-point results for nonexpansive maps go back to Browder and Kirk \cite{browder1965,kirk1965}, with key convergence principles and extensions (averaged operators, Fejér/quasi-Fejér monotonicity, alternating projections) surveyed in Bauschke–Combettes \cite{bauschke-combettes}. Our single-map clarification (\cref{prop:eq-single}, \cref{cor:classical}) reduces the “event-indexed” phrasing to the classical $T^N$-contractive case. For \emph{varying} operators, event-indexed \emph{block} contractivity (\cref{thm:main,thm:hetero}) gives a simple yet sharp envelope; this aligns with analyses that count effective contractive steps in block/intermittent regimes (e.g., randomized/coordinate updates, operator-splitting with periodic contraction).

\paragraph{Switching/dwell-time and randomized blocks (contrast).}
In switching-system stability, dwell-time constraints impose minimum on-periods for stabilizing modes; our bounded-gap hypothesis plays a similar counting role but at the level of abstract nonexpansive compositions, independent of system matrices. Randomized block methods (e.g., randomized coordinate or block-coordinate schemes) achieve convergence by ensuring sufficiently frequent contractive (or strongly averaged) updates in expectation; our heterogeneous envelope (\cref{thm:hetero}) provides a worst-case deterministic analogue driven by realized contractive blocks rather than averages.

\section{Scope of Validity and Non-Claims}

\begin{table}[htbp]
\caption{Scope summary.}
\label{tab:scope}
\centering
\footnotesize
\renewcommand{\arraystretch}{1.1}
\begin{tabular}{p{0.25\textwidth} p{0.24\textwidth} p{0.24\textwidth} p{0.23\textwidth}}
\toprule
\textbf{Result} & \textbf{Holds under} & \textbf{Not claimed / outside scope} & \textbf{Typical sufficient conditions} \\
\midrule
Anchored implication (Def.~\ref{def:anchored}; Props.~\ref{prop:reduction}--\ref{prop:no-synonym}) &
Fixed anchor $P$; semantics \eqref{eq:anchored-val}; reduction in $\mathcal{C}(P)$. &
No uniform Boolean synonym; \textbf{no completeness theorem}; no full residuation beyond $\mathcal{C}(P)$. &
Work inside $\mathcal{C}(P)$ (commutation $[E_B,P]=0$); Boolean reduction applies. \\
\addlinespace
Single-map events (Prop.~\ref{prop:eq-single}) &
$T$ nonexpansive; existence of a contractive power. &
No extra generality beyond $T^N$-contractive; after one event, all later powers are contractive. &
Show $\Lip(T^N)<1$ for some $N$ (e.g., strong monotonicity/averagedness). \\
\addlinespace
Varying-operator events (Thm.~\ref{thm:main}, \ref{thm:hetero}) &
Nonexpansive $T_t$; contractive blocks; $\mathrm{Fix}\neq\varnothing$. &
No converse; if no blocks occur, convergence may fail (Prop.~\ref{prop:noevents}). &
Projections $P_{C_t}$ with $\cap_t C_t\neq\varnothing$; shared minimizer; shared resolvent $J_{\gamma A_t}$ with common zero. \\
\addlinespace
Anchored convergence (Thm.~\ref{thm:anchored}) &
$TP=PT$; $(T|_{PH})^N$ contractive. &
No claim if $TP\neq PT$. &
Invariant subspace $PH$ and a contractive power of $T|_{PH}$. \\
\bottomrule
\end{tabular}
\end{table}

\section*{Authors' Synopsis of Key Claims (Alpay \& Kilictas)}
\begin{itemize}[leftmargin=1.4em,itemsep=0.35\baselineskip]
  \item \textbf{Heterogeneous blocks (factors $\boldsymbol{\lambda_k}$, gaps $\boldsymbol{N_k}$).} Product-form bounds and log-slope rates hold; see \cref{thm:hetero} and the examples in \cref{rem:hetero-examples}. Additional non-periodic schedules achieving the envelope exactly are provided above.
  \item \textbf{No uniform Boolean synonym.} No single Boolean $F(A,B)$ reproduces anchored implication uniformly over anchors; see \cref{prop:no-synonym} and the table therein.
  \item \textbf{Classical reduction in commuting contexts.} Inside $\mathcal{C}(P)$ (or $P=I$), anchored implication reduces to material implication with $E_{\neg A\lor B}=I-E_A+E_AE_B$; see \cref{prop:reduction}.
\end{itemize}

\section{Extensions and Future Work}\label{sec:future}
\textbf{Averaged operators.} Many algorithms use \emph{averaged} maps, $T=(1-\alpha)I+\alpha S$ with $S$ nonexpansive and $\alpha\in(0,1)$. Then $\Lip(T)\le 1$ and $T$ is firmly nonexpansive for $\alpha=\tfrac12$ in Hilbert spaces. A natural extension is \emph{block-averaged contractivity}: if each block is a composition/product of averaged maps whose aggregate is $\lambda$-contractive (e.g., via cocoercivity/strong monotonicity conditions), then \cref{thm:main,thm:hetero} apply with the same envelope. Further work could characterize minimal per-block conditions (e.g., strong averagedness in at least one step per window) that imply $\Lip(\Phi(n_k{:}n_{k-1}{+}1))<1$.

\section{Conclusion}

We clarified that for a single nonexpansive map, “event-indexed contractions” are equivalent to the classical $T^N$-contractive condition and thus do not introduce irregularity. The true gain comes with varying operators (switching/randomized), where contractive blocks can occur at irregular times and still yield uniqueness and explicit rates (provided a common fixed point exists). On the logic side, the anchored implication behaves classically in the commutant while providing a principled failure mode outside it. The envelope analysis extends to heterogeneous blocks with variable factors and gaps (\cref{thm:hetero}), and extending these ideas to averaged operators is a promising direction (\cref{sec:future}).

\appendix

\section{Appendix A: Reproducibility for Figs.~\ref{fig:fig1}--\ref{fig:fig2}}
Setup: $H=\mathbb{R}^2$, $x_0=(1,0)$, $z=0$, period $M=4$ with
$T_1=T_2=T_3=R_{\pi/2}$ and $T_4=\alpha I$ ($\alpha=0.8$), repeated.
Closed form:
\[
x_{4k+r}\ =\ R_{\pi/2}^{\,r}\,\big(\alpha^k R_{\pi/2}^{\,3k}x_0\big),\quad
\|x_{4k+r}-z\|=\alpha^k,\quad r\in\{0,1,2,3\}.
\]
Simple pseudocode (any language):
\begin{verbatim}
x = (1, 0)
alpha = 0.8
for n in 1..N:
    t = n % 4
    if t in {1,2,3}: x = R90(x)        # rotate by +90 degrees
    else:             x = alpha * x    # scalar contraction
    record norm(x)
\end{verbatim}

\section{Appendix B: Logic extensions}

\subsection*{B.1 Anchors as algebras}
One may replace the single anchor $P$ by a fixed von Neumann algebra $\mathcal{A}\subseteq B(H)$ and require the consequent to lie in its commutant: “$A \Rightarrow^{\mathrm{comm}}_{\mathcal{A}} B$” holds under the same semantics as \eqref{eq:anchored-val} but with the side condition $E_B\in\mathcal{A}'$. Inside $\mathcal{A}'$ the reduction to classical implication proceeds exactly as in \cref{prop:reduction}.

\subsection*{B.2 Effects (positive operators) via thresholds: mini-examples}
Let $A,B$ be \emph{effects}, i.e., $0\le A,B\le I$. For $0<\tau\le 1$, define threshold projections
\[
P_{A,\tau}:=\mathbf{1}_{[\tau,1]}(A),\qquad P_{B,\tau}:=\mathbf{1}_{[\tau,1]}(B).
\]
Define the \emph{$\tau$-anchored implication}
\[
v_\psi\!\big(A \Rightarrow^{\mathrm{comm},\tau}_{P} B\big)=1
\iff
\big(v_\psi(P_{A,\tau})=0\big)\ \text{or}\ \big(v_\psi(P_{A,\tau})=1,\ v_\psi(P_{B,\tau})=1,\ [B,P]=0\big).
\]

\paragraph{Mini-example 1 (commuting effects).}
Take $A=\mathrm{diag}(0.6,0.1)$, $B=\mathrm{diag}(0.7,0.2)$, $\tau=0.5$, and $P=\mathrm{diag}(1,0)$. Then
$P_{A,\tau}=\mathrm{diag}(1,0)$, $P_{B,\tau}=\mathrm{diag}(1,0)$, and $[A,P]=[B,P]=0$.
Hence
\[
E_{\,A \Rightarrow^{\mathrm{comm},\tau}_{P} B}\ =\ I - P_{A,\tau} + P_{A,\tau}P_{B,\tau}
=\mathrm{diag}(1,1)-\mathrm{diag}(1,0)+\mathrm{diag}(1,0)=I,
\]
matching the classical truth of $\neg A\lor B$ at this threshold.

\paragraph{Mini-example 2 (noncommuting effect; anchor failure).}
Let $U=\tfrac{1}{\sqrt{2}}\begin{psmallmatrix}1&1\\-1&1\end{psmallmatrix}$ and
$B=U\,\mathrm{diag}(0.9,0.1)\,U^\top$; set $A=\mathrm{diag}(1,0)$, $\tau=0.8$, $P=\mathrm{diag}(1,0)$.
Then $P_{A,\tau}=\mathrm{diag}(1,0)$, while $P_{B,\tau}$ is the rank-one projector onto $\mathrm{span}\{(1,-1)^\top\}$, so $[B,P]\neq 0$.
Even if a state $\psi$ satisfies $v_\psi(P_{A,\tau})=v_\psi(P_{B,\tau})=1$, the side condition fails and $A \Rightarrow^{\mathrm{comm},\tau}_{P} B$ evaluates to $0$; this mirrors the projection case.

\begin{proposition}[Reduction for commuting effects]
If $[A,P]=[B,P]=0$ (so $P_{A,\tau},P_{B,\tau}\in\mathcal{C}(P)$), then
\[
E_{\,A \Rightarrow^{\mathrm{comm},\tau}_{P} B}\ =\ I - P_{A,\tau} + P_{A,\tau}P_{B,\tau},
\]
i.e., the $\tau$-anchored implication reduces to the classical implication between the thresholded projections $P_{A,\tau}$ and $P_{B,\tau}$.
\end{proposition}

\end{document}